\documentclass[12pt,a4paper,leqno]{article}
\usepackage{amsmath} 
\usepackage{latexsym} 
\usepackage{amssymb}  

\newcommand{\irtb}{{\cal B}}

\newcommand{\irtg}{{\cal G}}

\newcommand{\irts}{{\cal S}}

\newcommand{\br}{{\bf R}}
\newcommand{\xk}{{X_\kappa}}
\newcommand{\rocfo}{{\varrho=cf(\varrho)}}
\newcommand{\roincalx}{{\varrho\in Cal(X)}}
\newcommand{\calx}{{Cal(X)}}

\newcommand{\cfalpha}{{cf(\alpha)}}

\newcommand{\upr}{\upharpoonright}
\newcommand{\hd}{{\hfill $\dashv$}}
\begin{document}
\title{Generic left-separated spaces and calibers}
\author{\\ I. Juh\'asz\ \footnote{\noindent Research   
supported by OTKA grant no. 37758.\newline
\indent 
{**}Research supported by The Israel 
Science Foundation. Publication 795.}\quad and \ S. Shelah $^{**}$}
\date{}
\maketitle

\begin{abstract}
In this paper we use a natural forcing to construct a 
left-separated topology on an arbitrary cardinal $\kappa$. 
The resulting left-separated space $\xk$ is also 
0-dimensional $T_2$, hereditarily Lindel\"of, and countably tight. 
Moreover if $\kappa$ is regular then $d(\xk)=\kappa$, hence 
$\kappa$ is not a caliber of $\xk$, while all other 
uncountable regular cardinals are. This implies that some 
results of [A] and [JSz] are, consistently, sharp.

We also prove it consistent that for every countable set $A$ 
of uncountable regular cardinals there is a hereditarily Lindel\"of 
$T_3$ space $X$ such that $\rocfo >\omega$ is a caliber of $X$ 
exactly if $\varrho\not\in A$. 
\end{abstract}



\setlength{\baselineskip}{18pt}
\parskip=10pt plus 5pt

\section*{ \S 1. Introduction} 

Let us start by recalling that a regular cardinal $\varrho$ 
is said to be a caliber of a topological space $X$ (in symbols: 
$\roincalx$) if among any $\varrho$ open subsets of $X$ there are 
always $\varrho$ many with non-empty intersection. Note that 
in this paper we restrict the notion of caliber to regular 
cardinals, although the definition does make sense for singular 
cardinals as well. Note also that $\roincalx$ implies that $X$ 
has no cellular family of size $\varrho$. Hence, as any infinite $T_2$ 
space has an infinite cellular family, for all spaces of 
interest we have 
$\calx\subset\br$, where $\br$ denotes the class of all 
uncountble regular cardinals. 

It is trivial to see that if $\rocfo >d(X)$ then $\roincalx$, 
moreover \v Sanin proved in [\v Sn] that, for any fixed $\varrho$, 
the property of 
spaces $\roincalx$
is fully productive. Consequently, for any cardinal $\kappa$ we have 
$Cal(2^\kappa)=Cal([0,1]^\kappa)=\br$, showing that the converse 
of the above relation between density and calibers is not valid. More 
precisely, no bound for the density of $X$ can be deduced 
from the fact that $X$ satisfies the condition
$\calx=\br$ that we also call \v Shanin's condition,
even for very nice (e. g. compact Hausdorff) spaces $X$. 

Such a converse, however, is valid if $X$ is a compact 
$T_2$ space of countable tightness, as was shown by 
\v Sapirovski\u\i\  in [\u Sap], see also [J1], 3.25. 
Indeed, in this case $\roincalx$ implies 
$d(X)<\varrho$ or, equivalently,
\[
\calx =[d(X)^+,\infty)
\]
where the interval on the right-hand side (just like in PCF 
theory) denotes an interval of regular cardinals.

More recently, in [A], Archangelski\u\i\ proved that if $X$ is 
Lindel\"of $T_3$ and countably tight and $\omega_1\in\calx$ then 
$d(X)\le 2^\omega$.
In [JSz] both \v Sapirovski\u\i\,'s and Archangelski\u\i\,'s 
results were strengthened and generalized, moreover, under CH, 
in the second result the conclusion $d(X)\le 2^\omega=\omega_1$ 
was improved to $d(X)=\omega$. Of course, this immediately led us 
to the question if the use of CH here is essential.

In the present note we give an affirmative answer to this 
question, in fact we show that Archangelski\u\i\,'s result is sharp 
for arbitrarily large values of the continuum $2^\omega$, 
even for hereditarily Lindel\"of (in short HL) $T_3$ spaces of 
countable tightness. The examples showing this will be 
obtained by forcing generic left-separated 
0-dimensional spaces in a natural way. Our methods will then be 
used to also solve some other problems raised in [JSz]. Moreover, 
we shall also prove the consistency of the statement that for any countable
subset $A$ of $\br$ there is a countably tight HL $T_3$ space 
$X$ such that
\[\calx =\br\setminus A\ .
\] 
This is in sharp contrast with the compact case. 

We do not know if there are similar consistency results for 
uncountable $A\subset\br$ and the following intriguing question also 
remains open: Is it provable in ZFC that a countably tight 
(hereditarily) Lindel\"of $T_3$ space $X$ satisfying \v Sanin's condition 
$\calx=\br$ is separable? 

Our notation and terminology follows [E] and [J1] in topology and 
[K] in forcing. 

\bigskip

\section*{ \S 2. Generic left-separated spaces} 

Let $\nu$ be an arbitrary limit ordinal and consider the suborder 
$P_\nu$ of the Cohen order $Fn(\nu^2,2)$ that consists of those 
$p\in Fn(\nu^2,2)$ which satisfy conditions (i) and (ii) below:
\newline (i) \ if $\langle\alpha,\alpha\rangle\in D(p)$ then 
$p(\alpha,\alpha)=1$;\newline
(ii) \ if $p(\alpha,\beta)=1$ then $\alpha\le\beta$.

Clearly, $P_\nu$ is a complete suborder of $Fn(\nu^2,2)$, hence it is 
CCC and thus preserves cardinals and cofinalities.

It is straight forward to check that for any pair $\langle\alpha,\beta
\rangle\in\nu^2$ the set 
\[
D_{\alpha,\beta}=\left\{p\in P_\nu:
  \langle\alpha,\beta\rangle\in D(p)\right\}
\]
is dense in $P_\nu$, consequently if $G \subset P_\nu$ is 
$P_\nu$-generic over $V$ then $$F=
\cup G :\nu^2\to 2,$$ i.e. $F$ defines a directed graph on $\nu$ by 
$F(\alpha,\beta)=1$ meaning that an edge goes from $\alpha$ to 
$\beta$. 

Now, in $V[G ]$, for any $\alpha\in\nu$ and $i\in 2$ let 
\[
U_{\alpha,i}=\left\{\beta\in\nu:F(\alpha,\beta)=1\right\}, 
\]
and $\tau_G $ be the (0-dimensional) topology on $\nu$ 
generated by the subbase
\[
\irts_G=\left\{U_{\alpha,i}:\alpha\in\nu,\ i\in 2\right\}.
\]
In other words, $\tau_G $ is the graph topology on 
$\nu$ determined by the directed graph $F$ in the sense of [J2] or [J3].

For all $\alpha\in\nu$ the minimal element of $U_{\alpha}$ is 
$\alpha$ and this shows that 
$\tau_G $ is left-separated in its natural
well-ordering. This immediately implies that $\tau_G  $ is 
$T_2$ and thus, by 0-dimensionality, also $T_3$. 

All finite intersections of the elements of $\irts_G $ form a base 
$\irtb_G $ of $\tau_G $. A typical element of 
$\irtb_G $ is of the form
\[
[\varepsilon]=\cap\left\{U_{\alpha,\varepsilon(\alpha)}:
    \alpha\in D(\varepsilon)\right\},
\]
where $\varepsilon\in Fn(\nu,2)$.

All this was easy. Let us now turn to the less obvious 
properties of the topology $\tau_G $. 

{\bf 2.1. Lemma.} \ {\em $\tau_G $ is {\rm HL}}.

{\bf Proof.} \ Assume, indirectly, that $p\in P_\nu$ forces 
that $\langle[\Dot\varepsilon_i]:i\in\omega_1\rangle$ are 
right-separating neighbourhoods of the points 
$\langle\Dot x_i:i\in\omega_1\rangle$ in $\nu$, where WLOG 
we may assume that $i<j$ implies $\Dot x_i<\Dot x_j$. Then for every
$i\in\omega_1$ there are $p_i\in P_\nu$, $\xi_i\in\nu$, and 
$\eta_i\in Fn(\nu,2)$ such that $p_i\le p$ and 
$p_i\Vdash\Dot x_i=\xi_i$ and $\Dot\varepsilon_i=\eta_i$. We 
may also assume that $D(p_i)=a_i^2$ for some 
$a_i\in[\nu]^{<\omega}$, moreover $\xi_i\in a_i$ and $D(\eta_i)\subset a_i$. 
By a standard $\Delta$-system and counting argument we can find
$i,j\in\omega_1$ with $i<j$ such that\newline 
(a) \ $p_i\upr (a_i\cap a_j)^2=p_j\upr (a_i\cap a_j)^2$, 
i.e.\ $p_i$ and $p_j$ are compatible as functions;\newline
(b) \ $\eta_i\upr a_i\cap a_j=\eta_j\upr a_i\cap a_j$;
\newline
(c) \ $\xi_i\in a_i\setminus a_j$, $\xi_j\in a_j\setminus a_i$, and 
$\xi_i<\xi_j$. 

Let us then define $q:(a_i\cup a_j)^2\to 2$ in such a way that $q\supset 
p_i\cup p_j$, moreover
\[
q(\alpha,\xi_j)=p_i(\alpha,\xi_i)
\]
if $\alpha\in D(\eta_i)\setminus a_j$ and $\alpha\le\xi_i<\xi_j$, and 
$q(\alpha,\beta)=0$ for every other pair $\langle\alpha,\beta\rangle
\in (a_i\cup a_j)^2\setminus (a_i^2\cup a_j^2)$.  
It is easy to see that $q\in P_\nu$ because it satisfies 
(i) and (ii), moreover $q(\alpha,\xi_j)=p_i(\alpha,\xi_i)$ holds 
for every $\alpha\in D(\eta_i)$, consequently
\[
q\Vdash\Dot x_j=\xi_j\in[\eta_i]=[\Dot\varepsilon_i],
\]
contradicting that $p\Vdash \Dot x_j\not\in [\Dot\varepsilon_i]$. \hd

\bigskip

Next we show that $\tau_G$ is countably tight.

\bigskip

{\bf 2.2. Lemma.} \ {\em $\tau_G $ has countably tightness.}

\bigskip

{\bf Proof.} \ Let us assume that for a $P_\nu$-name $\Dot A$ and some 
ordinal $\xi$ we have a condition $p\in P_\nu$ which forces 
$\xi\in{\Dot A}'$, i.e.\ that $\xi$ is an accumulation point of 
$\Dot A$. Since $\tau_G $ is left separated, we may also assume 
that $p\Vdash\xi<\Dot A$, i.e.\ $p$ forces that every element of 
$\Dot A$ is bigger than $\xi$. It can also be assumed that 
$\langle\xi,\xi\rangle\in D(p)$. 

Let $\lambda$ be a large enough regular cardinal such that $H(\lambda)$
contains ``everything in sight'', e.g.\ $P_\nu$, $\Dot A\in H(\lambda)$, 
etc. Fix a countable elementary submodel $N$ of $\langle
 H(\lambda),\in\rangle$ such that $\nu,\xi,p,\Dot A\in N$. Clearly, we shall 
be done if we can prove the following claim. 

\bigskip

{\bf Claim.} \ $p\Vdash\xi\in (N\cap \Dot A)'$. 

\bigskip

To see this, consider any $\varepsilon \in Fn(\nu,2)$ and let 
$q\le p$ be an arbitrary extension of $p$ in $P_\nu$ such that $D(q)=a^2$ 
with $a\in[\nu]^{<\omega}$, $D(\varepsilon)\subset a$, and 
$q\Vdash\xi\in [\varepsilon]$, i.e.\ $q(\alpha,\xi)=\varepsilon 
(\alpha)$ for every $\alpha\in D(\varepsilon)$. 

Then $q_N=q\cap N\in N$ is an extension of $p$ hence $q_N\Vdash 
|[\varepsilon_N]\cap\Dot A|\ge\omega$, 
where $\varepsilon_N=\varepsilon\cap N$. 
But we also have $q_N\in N$, hence $N\prec H(\lambda)$ implies 
that there is an extension $r\le q_N$ with $r\in N$
and an ordinal $x\in N\setminus a$ such that $D(r)=b^2$, $x\in b$, and 
\[
r\Vdash x\in N\cap\Dot A\cap [\varepsilon _N].
\]
Clearly $q\upr (a\cap N)^2=r\upr(a\cap N)^2=q_N$ and 
$D(q)\cap D(r)=(a\cap N)^2$, hence $q$ and $r$ are compatible as 
functions. We can thus define $q^*\supset q\cup r$ with the following additional 
stipulations: $q^*(\alpha,x)=\varepsilon(\alpha)$ whenever 
$\alpha\in D(\varepsilon)\setminus N$. Note that neihter $q$ nor 
$r$ is defined for a pair $\langle\alpha,x\rangle$ of this form 
because $x\not\in a\cap N$ and $\alpha\not\in N$. Also, $q^*\in P_\nu$ 
because (i) holds trivially and (ii) holds because if 
$\alpha\in D(\varepsilon)\setminus N$ and $\varepsilon(\alpha)=1$ then 
by $q(\alpha,\xi)=\varepsilon(\alpha)=1$ we have 
$\alpha\le\xi<x$. Finally, if $\alpha\in D(\varepsilon)\cap N$ then we 
have  
$q^*(\alpha,x)=r(\alpha,x)=\varepsilon_N(\alpha)=\varepsilon(\alpha)$ 
because $r\Vdash x\in[\varepsilon_ N]$, consequently
\[
q^*\Vdash x\in N\cap\Dot A\cap[\varepsilon]\ne\emptyset.
\]
This completes the proof of the claim and thus of lemma 2.2.
\hd

Note that lemmas 2.1 and 2.2 immediately yield us that $\langle 
\nu,\tau_G \rangle$ is a countably tight
$L$ space if $\nu\ge\omega_1$.

Our next lemma is the main result about calibers of $\tau_G $. In fact, 
for some applications to be given later, we formulate a slightly 
stronger result about calibers of initial segments of $\nu$ as 
subspaces of $\langle\nu,\tau_G \rangle$. 
So for $\alpha\le\nu$ we let $X_\alpha$
denote the subspace of 
$\langle\nu,\tau_G \rangle$ on $\alpha$. Note that for any $\beta\in
\nu\setminus\alpha$ we have $U_{\beta,1}
\cap\alpha=\emptyset$, consequently for any $\varepsilon\in Fn(\nu,2)$ 
we have either $[\varepsilon]\cap\alpha=\emptyset$ or 
$[\varepsilon]\cap\alpha=[\varepsilon\upr\alpha]\cap\alpha$. Therefore 
the trace of the base $\irtb_G $ on $\alpha$ can be written as
\[
\irtb_G \upr\alpha=
\{[\varepsilon]\cap\alpha:\varepsilon\in Fn(\alpha,2)\}\cup
   \{\emptyset\}.
\]

\bigskip

\textbf{2.3. Lemma.} \ {\em If $\alpha\le\nu$ is any limit ordinal 
and $\varrho$ is an uncountable regular cardinal with 
$\varrho<cf(\alpha)$ then $\varrho\in Cal(X_\alpha)$. 
Moreover, we also have $d(X_\alpha)=cf(\alpha)$ and so 
$Cal(X_\alpha)=\br\setminus\{cf(\alpha)\}$.} 

\bigskip

{\bf Proof.} \ By the above remark, to see the first part it clearly
suffices to show that whenever $p\Vdash\{\Dot\varepsilon_i:i\in\varrho\}
\subset Fn(\alpha,2)$ then for some $\xi\in\alpha$ there is a 
$q\le p$ such that
\[
q\Vdash|\{i\in\varrho:\xi\in[\Dot\varepsilon_i]\}|=\varrho.
\]

To see this, first we find for each $i\in\varrho$ an 
$\eta_i\in Fn(\alpha,2)$ and an extension $p_i\le p$ such that 
$p_i\Vdash\Dot\varepsilon_i=\eta_i$. We may also assume that 
$D(p_i)=a_i^2$ for some $a_i\in[\nu]^{<\omega}$ and 
$D(\eta_i)\subset a_i$ for all $i\in\varrho$. But then 
\[
|\cup\{a_i:i\in\varrho\}|\le\varrho<cf(\alpha),
\]
hence (the trace of) this union is bounded in $\alpha$. 
Consequently, there is an 
ordinal $\xi<\alpha$ with $a_i\cap\alpha<\xi$ for all $i\in\varrho$. 
Now extend each $p_i$ to a condition $q_i\in P_\nu$ such that 
$q_i(\alpha,\xi)=\eta_i(\alpha)$ for all $\alpha\in D(\eta_i)$. This 
is clearly possible because
\[
D(\eta_i)\subset a_i\cap\alpha<\xi\ .
\]
Note that then $q_i\Vdash\xi\in[\eta_i]=[\Dot\varepsilon_i]$.

Since $P_\nu$ is CCC and $q_i\le p$ for all $i\in\varrho$, 
there is a condition $q\in P_\nu$ with $q\le p$ such that
\[
q\Vdash|\{i\in\varrho : q_i\in\Dot G \}|=\varrho\ ,
\]
hence clearly
\[
q\Vdash|\{i\in\varrho : \xi\in[\Dot\varepsilon]\}|=\varrho\ ,
\]
which was to be shown.

To see that $d(X_\alpha)=cf(\alpha)$ first note that 
$d(X_\alpha)\ge cf(\alpha)$ is trivial because $X_\alpha$ is 
left-separated in its natural ordering. On the other hand, if 
$S\subset\alpha$ is any cofinal subset of $\alpha$ in the ground model 
$V$ then $S$ will be dense in $X_\alpha$. Indeed, it is again 
sufficient to show that $S\cap[\varepsilon]\ne\emptyset$ for every 
$\varepsilon\in Fn(\alpha,2)$, and this follows by a straight forward 
density argument. Consequently we have $d(X_\alpha)\le cf(\alpha)$, hence 
$d(X_\alpha)=cf(\alpha)$. 

Now, if $\varrho\in\br$ and $\varrho>d(X_\alpha)=cf(\alpha)$ then 
$\roincalx$, trivially. Finally, $\cfalpha\notin\calx$ is again obvious 
because $X_\alpha$ is left-separated. \hd

It is immediate from the above lemmas that if $\kappa$ is regular 
and $\kappa^\omega=\kappa$ then, in $V^{P_\kappa}$,
we have $2^\omega=\kappa$ 
and the space $X_\kappa$ is HL, 0-dimensional $T_2$, countably tight with 
$d(X_\kappa)=\kappa=2^\omega$, and $Cal(X_\kappa)=\br\setminus 
\{\kappa\}$. In particular, this shows that 
 Archangelski\u\i\,'s result from [A] saying that a Lindel\"of 
$T_3$ space $X$ with $\omega_1\in\calx$ 
satisfies $d(X)\le 2^\omega$ (or the more general corollary 1.2 of 
[JSz] saying that for such a space $X$ with 
$\varrho^+\in\calx$ we have $d(X)\le 2^\varrho$) is, at least consistently,
sharp.

Clearly, in a Lindel\"of space of countable tightness every free sequence
is countable.
Consequently, if we also have $\kappa>\omega_\omega$ then the space 
$X_\kappa$ establishes in addition that from 
corollary 1.5 of [JSz] (saying that if $X$ is a countably tight 
$T_3$ space with no free sequence of length $\omega_\omega$ and satisfying 
$\{\omega_n:0<n<\omega\}\subset\calx$  
then $X$ is separable
provided that
$\omega_\omega$ is strong limit) the assumption that $\omega_\omega$ 
be strong limit cannot be omitted.

With a little extra work we can deduce from our lemmas the following 
result showing that we have, again consistently, much more freedom 
in prescribing $\calx$ for Lindel\"of (even HL) and countably tight 
$T_3$ spaces than in the case of compact spaces of such kind.

\bigskip

\textbf{2.4. Theorem.} \ {\em Let $\kappa$ be any cardinal. Then, in
$V^{P_\kappa}$, for every countable subset $A$ of $\br\cap\kappa$ 
there is a HL and countably tight 0-dimensional $T_2$, hence $T_3$, space $X$
such that $\calx=\br\setminus A$.}

\bigskip

{\bf Proof.} \ 
For any $\varrho\in A $ let $X_\varrho$ be the subspace 
$\langle\varrho, \tau_G \upr\varrho\rangle$ as in 2.3 and then let
\[
X=\oplus\{X_\varrho:\varrho\in A\}
\]
be the (disjoint) topological sum of these subspaces. Since 
$A$ is countable, it is obvious that $X$ is HL, countably 
tight, and 0-dimensional $T_2$. For any $\varrho\in A$ then $X_\varrho$ 
is a clopen subspace of $X$, hence, by 2.3, we have 
$\varrho\not\in Cal(X_\varrho)$, implying that $\varrho\not\in \calx$
as well. 
On the other hand, if $\lambda\in\br\setminus A$ and $\irtg$ is a 
family of open sets in $X$ with $|\irtg|=\lambda$ then, again by the 
countability of $A$, there is a $\varrho\in A$ such that
$|\{G\in\irtg:G\cap X_\varrho\ne\emptyset\}|=\lambda$, hence 
by lemma 2.3 we have 
$\lambda\in Cal(X_\varrho)$ which implies that also 
$\lambda\in Cal(X)$. \hd

A natural question that we could not answer is if a similar result 
could be proved for 
uncountable sets $A$ of regular (uncountable) cardinals. 
Finally, our methods leave open the following very natural and 
interesting question formulated below.

\bigskip

\textbf{2.5. Problem.} \ Is it provable in ZFC that a 
Lindel\"of $T_3$ space $X$ of countable tightness satisfying 
\v Sanin's condition $\calx=\br$ is separable?

\vspace{2cm}

\noindent I. Juh\'asz\newline
A. R\'enyi Institute of Mathematics,\newline
Hungarian Academy of Sciences, \newline
Budapest, P.O.B. 127, 1364 Hungary\newline 
{\em e-mail address:} {\texttt {juhasz@renyi.hu}}

\noindent S. Shelah\newline
Institute of Mathematics,\newline
The Hebrew University of Jerusalem,\newline
91904 Jerusalem,\newline
Israel\newline
{\em e-mail address:} {\texttt {shelah@math.huji.ac.il}}


\vspace{2cm}

\begin{thebibliography}{99999}
\bibitem[A]{A} A. V. Archangelski\u \i, Projective $\sigma$-compactness
$\omega_1$-caliber, and $C_p$-spaces, Top.\ Appl.\ {\bf 104} (2000), 
13--16.
\bibitem[E]{E} R. Engelking, General Topology, Heldermann Verlag, Berlin, 
1989.
 
\bibitem[J1]{J1} I. Juh\'asz, Cardinal Functions - Ten Years Later, 
Math.\ Center Tracts, Vol. 123, Amsterdam, 1980.
 
\bibitem[J2]{J2} I. Juh\'asz, Cardinal Functions II, In: Handbook 
of General Topology, K. Kunen and J. E. Vaughan ed.s, 63--109, North--Holland,
Amsterdam, 1984.
 
\bibitem[J3]{J3} I. Juh\'asz, Cardinal Functions, In: 
Recent Progress in General Topology, M. Hu\v sek and 
J. van Mill ed.s, 417--441, North--Holland, Amsterdam, 1992.

\bibitem[JSz]{JSz} I. Juh\'asz and Z. Szentmikl\'ossy, Calibers, free 
sequences, and denstiy, Top.\ Appl.\ {\bf 119} (2002), 315--324.

\bibitem[K]{K} K. Kunen, Set Theory, North--Holland, Amsterdam, 1980. 

\bibitem[\v Sap]{vSap} B. E. \v Sapirovski\u\i, 
On tightness, $\pi$-weight, 
and related concepts, U\v c.\ Zap.\ Riga Univ.\ {\bf 3} (1976), 
88-89. (in Russian).

\bibitem[\v Sn]{vSn} N. A. \v Sanin, On the product of topological spaces, 
Trudy Math.\ Inst.\ Steklova {\bf 24} (1948) (in Russian).
\end{thebibliography}
\end{document}